\documentclass[10pt]{article}
\usepackage{amsmath}
\usepackage{amssymb}
\usepackage{theorem}
\usepackage{graphicx} 
\input supp-pdf.tex 

\newtheorem{theorem}{Theorem}
\newtheorem{proposition}{Proposition}

\newtheorem{lemma}[theorem]{Lemma}
\newtheorem{example}[theorem]{Example}
\newtheorem{remark}[theorem]{Remark}
\newcommand{\aut}{\mathrm{Aut}}
\newcommand{\aaut}{\mathrm{AAut}}
\newcommand{\caut}{\mathrm{CAut}}
\newcommand{\orb}{\mathrm{Orb}}
\newcommand{\cay}{\mathrm{Cay}}
\newcommand{\ul}[1]{\underline{#1}}

\title{Sporadic Examples of Directed \\ Strongly Regular Graphs Obtained by \\
Computer Algebra Experimentation \\ (Extended version)}
\author{\v Stefan Gy\"urki \\
{\small Matej Bel University, Faculty of Natural Sciences} \\
{\small Tajovsk\'eho 40, 974 01 Bansk\'a Bystrica, Slovak Republic, e-mail: {\texttt gypista@gmail.com}}
\and Mikhail Klin \\
{\small Department of Mathematics, Ben-Gurion University of the Negev} \\
{\small 84105 Beer Sheva, Israel, e-mail: {\texttt klin@cs.bgu.ac.il}}
}

\begin{document}

\maketitle


\begin{abstract}
We report about the results of the application of modern
computer algebra tools for construction of directed strongly
regular graphs. The suggested techniques are based on the
investigation of non-commutative association schemes and
Cayley graphs over non-Abelian groups. We demonstrate
examples of directed strongly regular graphs for 28 different
parameter sets, for which the existence of a corresponding digraph has
not been known before.
\end{abstract}

\section{Introduction}

This project is devoted to the computer algebra experimentation in
the area of algebraic graph theory, the part of mathematics on the
edge between graph theory, linear algebra, and group theory. The
main objects of interest in algebraic graph theory are highly
symmetric graphs, where level of symmetry might be measured both
on group-theoretical and purely combinatorial levels. Two books
\cite{bi} and \cite{gr} reflect impressive progress in this part
of mathematics.

Nowadays computer algebra tools, and especially {\sf GAP} (Groups, Algorithms, Programming -
a System for Computational Discrete Algebra \cite{ga}), together with
a few of its share packages, become an inalienable part of modern graph
theory and combinatorics. A~significant portion of striking combinatorial
structures was discovered and analyzed with the aid of a computer.
The main subject of interest in the presented text are directed strongly
regular graphs (briefly DSRGs), a natural generalization of a classical
(in algebraic graph theory) concept from simple to directed graphs.
The concept of a DSRG was suggested and investigated by A.~Duval in \cite{du}.
For a while it remained unnoticed, however, during last 15 years this class of
structures is becoming more and more popular.

The initial concept of a strongly regular graph (briefly SRG) has
a number of relatively independent origins of interest in such
diverse areas like design of statistical experiments, finite
geometries, applied permutation groups, and also complexity theory
of algorithms. Indeed, it is well-known that SRGs are usually
regarded as most sophisticated structures for the problems of
isomorphism testing of graphs and determination of the
automorphism group of graphs. The main combinatorial invariant of
a SRG is its parameter set, in the sense of \cite{bh}. Typically,
classification of SRGs is arranged for each parameter set
separately. Similar situation is also observed for DSRGs. However,
these structures appear even more frequently. For example, while
there are 36 parameter sets for SRGs on up to 50 vertices, this
number is 225 for DSRGs. On the other hand, the central problems
of the identification of DSRGs and determination of their symmetry
are on the same level of difficulty as it appears for the
classical case of strongly regular graphs.

In this context, DSRGs provide, in comparison with SRGs much more wide training
polygons for the experts in the complexity theory which allows more diversity
(undirected versus directed) for investigated graphs.

The previous experience (earned, in particular, by M. Klin and his
coauthors) shows that a clever use of computers helps to discover
new examples of DSRGs and after that to reach an honest
theoretical generalization of the detected structures. This line
of activity stimulated the authors to join their efforts in a new
attempt. At this stage, we are concentrating on the association
schemes as possible origins of new DSRGs. Namely, we wish to
consider any association scheme $\mathcal M$, for which a suitable
union of classes provides a DSRG, preferably new, moreover, with a
new parameter set.

This text is an extended version of a paper (with the same title)
to appear soon in a special volume (CASC 2014) of LNCS series. 

The paper is organized as follows. In Section~2, the necessary basic notions
are introduced. In Section~3, we describe our approach
to the problem of finding new directed strongly regular graphs
using  computer algebra experimentation. In Section~4, the mentioned strategies
are explained with enough rigorous details and the results of different approaches
are reported.  In Section~5, the results of a classical strategy using
Cayley graphs are submitted.  We conclude with a discussion and summary of new
graphs, being discovered. Information, provided in the Appendix, allows one to reconstruct
(with the aid of a computer) all new DSRGs discovered by us. 

\section{Preliminaries}

Below we present brief account of most significant concepts
exploited in the paper. We refer to \cite{bh} and \cite{pw} for
more information.

\subsection{General Concepts}

A \emph{simple graph} $\Gamma$ is a pair $(V,E)$, where $V$ is a finite set of \emph{
vertices}, and $E$ is a set of 2-subsets of $V$ which are called \emph{edges}.

A \emph{directed graph} (briefly digraph) $\Gamma$ is a pair $(V,R)$ where $V$ is the
set of \emph{vertices} and $R$ is a binary relation on $V$, that is a subset of the set $V^2$
of all ordered pairs of elements in $V$. The pairs in $R$ are called \emph{directed arcs} or \emph{darts}.
The vertex set of $\Gamma$ is denoted by $V(\Gamma)$ and the dart set is denoted by $R(\Gamma)$.

A \emph{balanced incomplete block design} (BIBD) is a pair $(\mathcal P,\mathcal B)$
 where $\mathcal P$ is the point set of cardinality $v$,
 and $\mathcal B$ is a collection of $b$ $k$-subsets of $\mathcal P$
(\emph{blocks}) such that each element of $\mathcal P$ is
contained in exactly $r$ blocks and any 2-subset of $V$ is
contained in exactly $\lambda$ blocks. The numbers $v, b, r, k$,
and $\lambda$ are \emph{parameters} of the BIBD. From the
parameters $v,k,\lambda$ the remaining two are determined
uniquely, therefore, we use just the triplet of parameters
$(v,k,\lambda)$ for a BIBD.

For any finite group $H$, the {\it group ring} $\mathbb ZH$ is
defined as the set of all formal sums of elements of $H$, with
coefficients from $\mathbb Z$. Let $X$ denote a non-empty subset
of $H$. The element $\sum_{x\in X}x$ in $\mathbb ZH$ is called a
{\it simple quantity}, and it is denoted as $\ul{X}$. Suppose now
that $e\notin X$, where $e$ is the identity element of the group
$H$. Then the digraph $\Gamma=\cay(H,X)$ with vertex set $H$ and
dart set $\{(x,y): x,y\in H, yx^{-1}\in X\}$ is called the {\it
Cayley digraph over $H$ with respect to~$X$.}

\subsection{Strongly Regular Graphs}

A graph $\Gamma$ with adjacency matrix $A=A(\Gamma)$ is called
\emph{regular}, if there exists a positive integer $k$ such that
$AJ=JA=kJ$, where $J$ is the all-one matrix. The number $k$ is
called \emph{valency} of $\Gamma$. A simple regular graph with
valency $k$ is said to be \emph{strongly regular} (SRG, for short)
if there exist integers $\lambda$ and $\mu$ such that for each
edge $\{u,v\}$ the number of common neighbors of $u$ and $v$ is
exactly $\lambda$; while for each non-edge $\{u,v\}$ the number of
common neighbors of $u$ and $v$ is equal to $\mu$. Previous
condition can be rewritten equivalently into the equation
$A^2=kI+\lambda A+\mu(J-I-A)$ using the adjacency matrix of
$\Gamma$. The quadruple $(n,k,\lambda,\mu)$ is called the
\emph{parameter set} of an SRG $\Gamma$.

\subsection{Directed Strongly Regular Graphs}

A possible generalization of the notion of SRGs for directed graphs was given by Duval~\cite{du}.
While the family of SRGs has been well-studied in the algebraic graph theory cf. \cite{bh},
the directed version has not received enough attention.

A \emph{directed strongly regular graph} (DSRG) with parameters $(n,k,t,\lambda,\mu)$ is a
regular directed graph on $n$ vertices with valency $k$, such that every vertex is incident with~$t$
undirected edges, and the number of paths of length 2 directed from a vertex $x$ to another vertex $y$ is
$\lambda$, if there is an arc from $x$ to $y$, and $\mu$ otherwise.
In particular, a DSRG with $t=k$ is an SRG, and a DSRG with $t=0$ is a doubly regular tournament.
Throughout the paper we consider only DSRGs satisfying $0<t<k$, which are called \emph{genuine}
DSRGs.

The adjacency matrix $A=A(\Gamma)$ of a DSRG with parameters $(n,k,t,\lambda,\mu)$,
 satisfies $AJ=JA=kJ$ and $A^2=tI+\lambda A+\mu(J-I-A)$.

\medskip

\begin{figure}[htb]
\begin{center}
\includegraphics{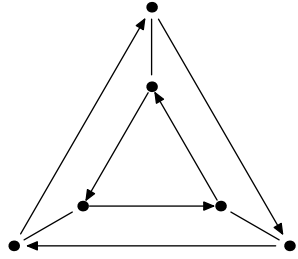}
\end{center}
\caption{The smallest genuine DSRG.}
\label{fig1}
\end{figure}

\noindent{\bf Example 1.}
The smallest example of a DSRG is appearing on 6 vertices. Its
parameter set is $(6,2,1,0,1)$ and it is depicted in
Fig.~\ref{fig1}.

\medskip

\begin{remark}
In this paper, we are using for DSRG's 5-tuple of parameters in
the order $(n,k,t,\lambda,\mu)$, however, in several other papers
the order $(n,k,\mu,\lambda,t)$ is used.
\end{remark}

\begin{proposition}[\cite{du}]
If $\Gamma$ is a DSRG with parameter set $(n,k,t,\lambda,\mu)$ and adjacency matrix $A$,
then the complementary graph $\bar\Gamma$ is a DSRG with parameter set $(n,\bar k,\bar t,
\bar\lambda,\bar\mu)$ with adjacency matrix $\bar A=J-I-A$, where
\begin{eqnarray*}
\bar k &=& n-k+1 \\
\bar t &=& n-2k+t-1 \\
\bar \lambda & =& n-2k+\mu-2 \\
\bar \mu &=& n-2k+\lambda.
\end{eqnarray*}
\label{prop2}
\end{proposition}

\begin{remark}
Proposition \ref{prop2} allows us to restrict our search for the DSRGs with $2k<n$,
due to complementation, and clearly a discovery of a DSRG with new parameter set
implies a discovery of a DSRG on the complementary parameter set.
As a consequence, throughout the paper we display just the parameter sets satisfying $2k<n$.
\end{remark}

For a directed graph $\Gamma$ let $\Gamma^T$ denote the digraph obtained by reversing
all the darts in $\Gamma$. Then $\Gamma^T$ is called the \emph{reverse} of $\Gamma$. In other words,
if $A$ is the adjacency matrix of~$\Gamma$, then $A^T$ is the adjacency matrix of $\Gamma^T$.

The following proposition was observed by Ch. Pech, and presented in \cite{km}:

\medskip

\begin{proposition}[\cite{km}]
Let $\Gamma$ be a DSRG. Then the graph $\Gamma^T$ is also a DSRG with the same parameter set.
\label{prop1}
\end{proposition}

We say that two DSRGs $\Gamma_1$ and $\Gamma_2$ are
\emph{equivalent}, if $\Gamma_1\cong\Gamma_2$, or
$\Gamma_1\cong\Gamma_2^T$, or $\Gamma_1\cong\bar\Gamma_2$, or
$\Gamma_1\cong\bar\Gamma_2^T$; otherwise they are called
\emph{non-equivalent}. (In other words, $\Gamma_1$ is equivalent
to $\Gamma_2$ if and only if $\Gamma_1$ is isomorphic to
$\Gamma_2$ or to a graph obtained from $\Gamma_2$ via reverse and
complementation.) From our point of view the interesting DSRGs are
those which are non-equivalent.

The parameters $n,k,t,\lambda,\mu$ are not independent. Relations
to be satisfied for such parameter sets are usually called
\emph{feasibility conditions}. Most important and, in a sense,
basic conditions are the following (for their proof see
\cite{du}):

\begin{equation}
k(k+\mu-\lambda)=t+(n-1)\mu. \label{c1}
\end{equation}

There exists a positive integer $d$ such that:
\begin{eqnarray}
d^2&=&(\mu-\lambda)^2+4(t-\mu) \label{c2}\\
d &\mid& (2k-(\mu-\lambda)(n-1)) \label{c3}\\
n-1 &\equiv& \frac{2k-(\mu-\lambda)(n-1)}{d}  \pmod 2 \label{c4}\\
n-1 &\geq &\left|\frac{2k-(\mu-\lambda)(n-1)}{d}\right|.\label{c5}
\end{eqnarray}

Further:
\[
\begin{array}{rcccl}
0 & \leq & \lambda < t & < & k \\
0 & < & \mu \leq t & < & k \\
-2(k-t-1) & \leq & \mu  -  \lambda & \leq & 2(k-t).
\end{array}
\]

We have to mention that for a feasible parameter set
it is not guaranteed that a DSRG with that parameter set does exist.
A feasible parameter set for which at least one DSRG $\Gamma$ exists is called \emph{realizable},
otherwise \emph{non-realizable}.
The smallest example of a non-realizable parameter set is $(14,5,4,1,2)$,
what was shown in \cite{km}.

\subsection{Coherent Configurations and Association Schemes}

Under a \emph{color graph} $\Gamma$ we mean an ordered pair
$(V,\mathcal R)$, where $V$ is a set of vertices and $\mathcal R$
a partition of $V\times V$ into binary relations. The elements of
$\mathcal R$ are called \emph{colors}, and the number of
colors is the \emph{rank} of $\Gamma$. In other words, a color graph
is an edge-colored complete directed graph with loops, whose arcs are colored
by the same color if and only if they belong to the same binary relation.

A \emph{coherent configuration} is a color graph $\mathcal M=(\Omega,
\mathcal R)$, $\mathcal R=\{R_i \mid i\in I\}$, such that the following
axioms are satisfied:
\begin{itemize}
\item[(i)] The diagonal relation $\Delta_{\Omega}=\{(x,x)\mid x\in\Omega\}$
is a union of relations $\cup_{i\in I'}R_i$, for a suitable subset $I'
\subseteq I$.
\item[(ii)] For each $i\in I$ there exists $i'\in I$ such that
$R_i^T=R_{i'}$, where $R_i^T=\{(y,x)\mid (x,y)\in R_i\}$ is the relation
transposed to $R_i$.
\item[(iii)] For any $i,j,k\in I$, the number $p_{i,j}^k$ of elements $z\in
\Omega$ such that $(x,z)\in R_i$ and $(z,y)\in R_j$ is a constant
depending only on $i,j,k$, and independent of the choice of
$(x,y)\in R_k$.
\end{itemize}
The numbers $p_{i,j}^k$ are called \emph{intersection numbers}, or
sometimes \emph{structure constants} of $\mathcal M$.
A coherent configuration $\mathcal M$ is called \emph{commutative},
if for all $i,j,k\in I$ we have $p_{ij}^k=p_{ji}^k$, otherwise \emph{non-commutative}.

Let $(G,\Omega)$ be a permutation group. $G$ acts naturally on $\Omega\times\Omega$ by
$(x,y)^g=(x^g,y^g)$. The orbits of this action are called \emph{2-orbits} (or \emph{orbitals})
of $(G,\Omega)$, and denoted by $2{-}\orb(G,\Omega)$.
It is easy to check that $(\Omega,2{-}\orb(G,\Omega))$ is a coherent configuration
for every permutation group $(G,\Omega)$. The coherent configurations which appear
in this manner are called \emph{Schurian}, otherwise \emph{non-Schurian}.

An \emph{association scheme} $\mathcal M=(\Omega,\mathcal R)$ is a
\emph{homogeneous} coherent configuration, i.e., where the
diagonal relation $\Delta_{\Omega}$ does belong to $\mathcal R$.
Hence, a very important source of association schemes are
transitive permutation groups, since their 2-orbits form a
homogeneous coherent configuration.

Let $\mathcal M$ be a coherent configuration of rank $r$.
To each relation $R_i$ in $\mathcal M$ we can assign a
0--1-matrix $A_i$ such that $(A_i)_{xy}=1 \iff (x,y)\in R_i$.
Then clearly $\sum_{i=1}^r A_i=J$ and $A_iA_j=\sum_{k=1}^r p_{ij}^kA_k$.
Matrices $A_1,\ldots,A_r$ generate an algebra $\mathcal W$ over $\mathbb C$,
which is called \emph{coherent algebra} of rank $r$ and degree $n$, and we write
$\mathcal W=\langle A_1,\ldots,A_r\rangle$.

Let $\mathcal W_1,\mathcal W_2$ be two coherent algebras of order
$n$. Then $\mathcal W_1\cap \mathcal W_2$ is again a coherent
algebra, therefore, there exists a unique minimal coherent algebra
$\mathcal W$ containing a given set $\{M_1,\ldots,M_t\}$ of
0--1-matrices of order $n\times n$. This algebra is called
\emph{coherent closure} of $M_1,\ldots,M_t$ and it is denoted
$\langle\langle M_1,\ldots,M_t\rangle\rangle$. In particular, to a
DSRG $\Gamma$ with adjacency matrix $A$, we can associate the
coherent closure $\mathcal W(\Gamma)=\langle\langle
A\rangle\rangle$.

To each coherent configuration $\mathcal M$, we can assign three
groups: $\aut(\mathcal M)$, $\caut(\mathcal M)$ and
$\aaut(\mathcal M)$. The (combinatorial) \emph{group of
automorphisms} $\aut(\mathcal M)$ consists of the permutations
$\phi:\Omega\to\Omega$ which preserve the relations, i.e.,
$R_i^{\phi}=R_i$ for all $R_i\in\mathcal R$. The \emph{color
automorphisms} preserve relations setwise, i.e., for
$\phi:\Omega\to\Omega$ we have $\phi\in\caut(\mathcal M)$ if and
only if for all $i\in I$ there exists $j\in I$ such that
$R_i^{\phi}=R_j$. An \emph{algebraic automorphism} is a bijection
$\psi:\mathcal R\to \mathcal R$ which satisfies
$p_{ij}^k=p_{i^{\psi}j^{\psi}}^{k^{\psi}}$. We refer to \cite{pw}
for a discussion of these concepts.

Graphs and digraphs can be regarded as binary relations, while
association schemes are collections of binary relations in the
sense of our definition. Therefore, it is natural to ask:

\medskip

\noindent{\bf Question 1.}
Assume $\mathcal M$ is an association scheme of order $n$.
Can we obtain a DSRG on $n$ vertices as a union of suitable
classes in $\mathcal M$?

\medskip

It turns out that there is no standard easy way to reply to the
question.

A very important necessary condition posed for the initial association scheme $\mathcal M$
was given in \cite{km}:

\begin{theorem}[\cite{km}]
Let $\Gamma$ be a genuine directed strongly regular graph. Then
the coherent closure $\mathcal W(\Gamma)$ is non-commutative, and
its rank is at least $6$.
\end{theorem}

In other words, we have to consider non-commutative association schemes
of rank at least 6, when we are searching for directed strongly regular graphs
as unions of relations in a prescribed association scheme.

Significant part of our results was achieved following the strategy of
creating suitable non-commutative association schemes and taking
unions of their relations.

\section{General Approach to the Computer Experimentation}

\subsection{Main Methodology}
Assume that $\mathcal M=(\Omega, \mathcal R)$ is an association
scheme and $G=\aut(\mathcal M)$. Let $r$ be the rank of $\mathcal
M$, thus, $\mathcal M$ has $r-1$ classes. In many cases below, $G$
acts transitively on $\Omega$. Moreover, $\mathcal M$ is the
Schurian scheme obtained from permutation group $(G,\Omega)$,
however, this restriction is not obligatory in the framework of
the described approach.

Let $\Gamma$ be a putative DSRG (with order $n$), which is
obtained via union of suitable classes of $\mathcal M$. Then clearly
one has to inspect $2^{r-1}$ possible unions.

First evident restriction is to look simultaneously for all possible
parameter sets of DSRGs of order $n$; recall that this data is
available at \cite{ab}. Typically, in this  project,
our attention was restricted only to the open parameter sets.

At the second step, one has to consider multisets of valencies of
symmetric and antisymmetric classes in $\mathcal M$ and to find in
advance which subsets of classes of $\mathcal M$ may in principle
provide a mixed graph with prescribed pair of valencies $(t,
k-t)$, respectively. Getting such a list is a simple case of the
famous knapsack problem, however, we were using a very naive
approach to provide all solutions.

In order to eliminate in further search duplicates of isomorphic
graphs it might be helpful to work with the representatives of
orbits on sets of relations in the action of the color group
$\caut(\mathcal M)$. Sometimes, preliminary sorting with the aid
of the action of the algebraic group $\aaut(\mathcal M)$ might be
also of help. Nevertheless, according to the gathered practical
experience, in most of the cases the order of $\aaut(\mathcal M)$
is relatively small, thus, we are facing exponential complexity in
cases of schemes with relatively small valencies. Therefore, we
decided to restrict our systematic attempts just to the
association schemes of rank not larger than 25.

Finally, for each selected ``suspicious'' union of relations from $\mathcal M$
we have to check whether it is providing a DSRG or not. Here use of the known
structure constants of $\mathcal M$ is very crucial: indeed, instead of inspection
of the adjacency matrix $A(\Gamma)$ of a putative graph $\Gamma$ we arrange calculations
with the tensor of structure constants of $\mathcal M$.

The computational scheme outlined above is, in a sense, the ideal
plan of activities, which were arranged in the course of
computations. In many cases, we preferred to use ingredients of a
brute-force approach, rather than to being involved in a more
sophisticated programming. Since in many cases it was impossible
to execute an exhaustive search, it was substituted by an ad hoc
selection of simple ``promising'' subsets of candidates.

\subsection{Computer Tools}

We run all computations in the software {\sf GAP} \cite{ga} with its
share packages {\sf GRAPE} \cite{lh} together with {\sf nauty}
\cite{na} for computation with graphs;
an unpublished package {\sf COCO-II} \cite{sv} written by S. Reichard
for computations with association schemes and coherent configurations;
and the package {\sf SetOrbit} \cite{so} written by Ch. Pech and S. Reichard,
and documented in \cite{pr}, for finding representatives of orbits of group actions on sets of
various size.

In addition, some ad hoc computational tricks were used from time to time, like
to exploit a simple variation of the calculation of the coherent closure of an
auxiliary graph, which is related to the putative DSRG $\Gamma$, as well as some
helpful functions for the calculations with association schemes borrowed from the
site \cite{hm}.

\subsection{Sources for Association Schemes}

Recall that first open parameter set for a DSRG appears for order $n=22$.
With growing of $n$ the fraction of open parameter sets is becoming more
essential. This dictated our strategy in the selection of candidates for
association schemes being considered. In what follows, we report only about
successful attempts, resulted in discovery of graphs with open parameter sets.
However, as a byproduct, many graphs with known parameter sets were also considered
(their comparison with known ones remains as one of tasks for a more systematic
approach in the future).

Roughly speaking, we distinguish a few different typical origins in our search:
\begin{itemize}
\item use of existing catalogues of association schemes;
\item inspection of groups of automorphisms of some ``famous'' vertex-transitive
graphs;
\item consideration of incidence structures;
\item investigation of Cayley graphs.
\end{itemize}

In the next section, we are paying reasonable attention to a more
detailed discussion of each of these approaches.

\section{Unions of Relations in Association Schemes}

Here we consider several strategies for finding non-commutative
association schemes which serve as input for searching new DSRGs.

\subsection{Search Using Catalogue of Small Association Schemes}

For executing our strategy it is enough to consider
non-commutative association schemes of small order at the first stage.
They are systematically arranged according to their order and rank in the
catalogue of Hanaki and Miyamoto \cite{hm}.

The number of new parameter sets, for which we succeed, using
exactly this approach, is $12$, see Table \ref{tabS} in Summary.
For several parameter sets, we have found a few non-equivalent
DSRGs. Table \ref{tab1} contains just the digraphs which are
mutually non-equivalent. In this Table \ref{tab1}, we display
sufficient portion of information for reconstructing discovered
DSRGs using the catalogue of association schemes by Hanaki and
Miyamoto.

\begin{remark}
We noticed that in \cite{hm} ``class'' of association schemes is used instead
of their ``rank''. Clearly, the number of classes is less by one than the rank.
\end{remark}

From Table \ref{tab1} it is easy to observe that graph nr. 13 is a (spanning)
subgraph of graphs nr. 25 and 26; while nr. 14 is a subgraph of nr. 19.

\begin{table}[h]
\caption{DSRGs from small association schemes (Used abbreviations: AS -- association scheme, nr.cat -- 
number of the AS in the catalogue)}
\label{tab1}
\begin{center}
\begin{tabular}{|r|l|l|c|c|c|}
\hline
nr. & $(n,k,t,\lambda,\mu)$ & Union of relations & AS order & AS rank & nr.cat \\ 
\hline
1 & $(30,13,11,6,5)$ & $1,2,4,6,8$ & 30 & 11 & 184 \\ 
2 & $(36,13,7,4,5)$ & $4,5,6,9$ & 36 & 11 & 49 \\ 
3 & $(36,13,7,4,5)$ & $1,2,3,4,8,12,14,16$ & 36 & 20 & 28 \\ 
4 & $(36,13,7,4,5)$ & $1,2,3,4,8,12,14,17$ & 36 & 20 & 28 \\ 
5 & $(36,13,7,4,5)$ & $1,2,3,4,5,10,13,14$ & 36 & 20 & 30 \\ 
6 & $(36,13,7,4,5)$ & $1,2,3,4,6,8,12,18$ & 36 & 20 & 40 \\ 
7 & $(36,13,7,4,5)$ & $1,2,3,4,6,8,13,16$ & 36 & 20 & 40 \\ 
8 & $(36,13,11,2,6)$ & $1,3,5,7,10$ & 36 & 13 & 57 \\ 
9 & $(45,16,8,5,6)$ & $1,3,5,8$ & 45 & 10 & 18 \\ 
10 & $(45,16,8,5,6)$ & $2,3,5,8$ & 45 & 10 & 18 \\ 
11 & $(50,16,10,3,6)$ & $1,2,6,9,13$ & 50 & 14 & 9 \\ 
12 & $(50,23,13,10,11)$ & $1,2,6,8,10,12$ & 50 & 14 & 17 \\ 
13 & $(54,8,3,2,1)$ & $1,6,8,12$ & 54 & 18 & 103 \\ 
14 & $(54,16,12,6,4)$ & $4,5,10,12,16$ & 54 & 18 & 109 \\ 
15 & $(54,19,9,6,7)$ & $1,3,4,6,11,14,16$ & 54 & 18 & 111 \\ 
16 & $(54,19,9,6,7)$ & $1,3,4,6,11,14,17$ & 54 & 18 & 111 \\ 
17 & $(54,20,16,6,8)$ & $2,3,4,5,11,13,14$ & 54 & 18 & 109 \\ 
18 & $(54,20,16,6,8)$ & $8,9,11,13,14$ & 54 & 18 & 109 \\ 
19 & $(54,21,17,8,8)$ & $1,2,3,4,5,10,12,16$ & 54 & 18 & 109 \\ 
20 & $(54,21,17,8,8)$ & $1,8,9,10,12,16$ & 54 & 18 & 109 \\ 
21 & $(54,25,14,11,12)$ & $1,2,4,6,10,12,14,16$ & 54 & 18 & 106 \\ 
22 & $(54,25,14,11,12)$ & $1,2,4,6,10,12,15,17$ & 54 & 18 & 106 \\ 
23 & $(54,25,14,11,12)$ & $1,2,4,6,10,13,15,16$ & 54 & 18 & 106 \\ 
24 & $(54,25,14,11,12)$ & $1,2,4,6,11,13,15,17$ & 54 & 18 & 106 \\ 
25 & $(54,25,14,11,12)$ & $1,2,6,8,10,12,14,16$ & 54 & 18 & 103 \\ 
26 & $(54,25,14,11,12)$ & $1,2,6,8,10,12,15,17$ & 54 & 18 & 103 \\ 
27 & $(54,25,14,11,12)$ & $1,2,6,8,10,13,14,17$ & 54 & 18 & 103 \\ 
28 & $(54,25,14,11,12)$ & $1,2,6,8,11,13,14,16$ & 54 & 18 & 103 \\ 
\hline
\end{tabular}
\end{center}
\end{table}

\subsection{Actions of Group of Automorphisms of Graphs}

J\o rgensen in \cite{lk} and \cite{jr} announced the existence of
a DSRG with the parameter set $(108,10,3,0,1)$. The author
provided us the adjacency matrix of this new digraph, and we
managed to explain it in terms of unions of relations in the
Schurian association scheme of the group of automorphism of the
Pappus graph in the action on the ordered triples of its vertices.
For more details see \cite{gk}.

This successful attempt inspired us to go ahead in a similar
spirit. In fact, we investigated actions of the group of
automorphisms of several symmetric graphs on certain orbits of
various $k$-sets and $k$-tuples. Usually, due to high time and
space complexity, we took just $k\in\{2,3,4\}$. Restricting group
action to an orbit we ensure that the resulted action is
transitive on it, and from this action we create the Schurian
association scheme. When it passes the test for being
non-commutative, then there is sense to execute the search for
DSRGs as unions of relations in these schemes. Once more, due to
high time-complexity, we restricted ourselves just for the cases
when the rank was not greater than $25$ and the size of the orbit
not greater than $110$. Therefore, our search is far from being
exhaustive. If one goes higher with the rank, then he could
probably find new DSRGs.

One can find origins of this strategy in \cite{fp}, Example 3.4.
The authors took the lattice graph on 9 points and investigated an
action of a subgroup of its group of automorphism on the edges.
Our strategy is a slight generalization of it, since we do not
consider only the pairs of two adjacent vertices, but also actions
on any 2-sets, 3-sets, ordered pairs, ordered triplets of vertices
and sometimes on 4-sets, ordered quadruples.

Using our strategy, we succeed in the following cases (the starting ``famous'' graphs are
available via Internet, e.g. from the home page of A. Brouwer):

\begin{itemize}
\item from the Petersen graph we obtained a DSRG$(60,13,5,2,3)$;
\item from the Shrikhande graph we get a DSRG$(48,10,6,2,2)$ and $(48,13,7,2,4)$;
\item from the Heawood graph we obtained DSRG $(84,31,17,12,11)$,
$(84,29,19,6,12)$ and $(84,39,27,18,18)$, see also Section 4.3;
\item from the unique SRG$(21,10,3,6)$ we obtained a DSRG$(105,36,16,11,13)$;
\end{itemize}

Explicit descriptions of these digraphs are shown in Appendix.

\subsection{Actions of Group of Automorphisms of Combinatorial Designs}

Let us now start from a block design $\mathcal D=(\mathcal
P,\mathcal B)$ with the point set $\mathcal P$ and block set
$\mathcal B$, let $I=I(\mathcal D)$ be the Levi graph of $\mathcal
D$, that is the graph with vertex set $\mathcal P\cup\mathcal B$
and two vertices being adjacent if and only if the corresponding
elements of $\mathcal D$ are incident. Clearly, $I(\mathcal D)$ is
a bipartite graph. The group $\aut(I(\mathcal D))$ either
coincides with the group $\aut(\mathcal D)$, or it is twice larger
(the latter corresponds to the case when $\mathcal D$ is a
symmetric self-dual design).

For a number of designs $\mathcal D$, we investigated the action
of the group $G=\aut(I(\mathcal D))$ on certain orbits of various
$k$-subsets and $k$-tuples of vertices of $I$. The same
limitations for values of $k$, order and rank of related Schurian
association scheme (like in previous section) remain valid. The
execution of search for DSRGs has been started, provided the
appearing Schurian schemes were non-commutative.

Using this strategy we succeed in the following cases:

\begin{itemize}
\item Considering the unique $(7,3,1)$-design is equivalent of consideration of the Heawood
graph in the previous section, since the Levi graph of the
$(7,3,1)$-design is the Heawood graph;

\item from the unique $(9,3,1)$-design we get DSRGs with parameter sets $(72,20,14,4,6)$,
$(72,21,15,6,6)$, $(72,22,9,6,7)$;
\item from the $(10,4,2)$-design with group of automorphisms of order $720$ we get
DSRGs with parameter sets $(60,26,20,10,12)$ and $(90,28,16,10,8)$;
\item from the $(15,3,1)$-design which has group of automorphisms of order $288$ we get
DSRG with parameter set
$(72,26,10,8,10)$;
\item from the $(15,7,3)$-design with group of automorphisms of order $1152$ we get
a DSRG with parameter set $(72,19,11,2,6)$;
\item from the $(25,4,1)$-design with group of automorphisms of order $504$ we get
a DSRG with parameter set $(63,22,10,7,8)$.
\end{itemize}

In all these cases, we refer to the description of block designs
provided in~\cite{cd}.

Explicit descriptions of all digraphs constructed in this subsection are shown in Appendix.

\begin{example}
Consider the unique $(7,3,1)$-design $F$, that is the Fano plane.
In this case, the group $G$ of order 336 acts transitively on the
vertex set of graph $I(F)$ of size 14. Let us consider
configuration, which consists of two lines, their intersection
point and another point in one line not belonging to the other
line. Clearly, there are ${7\choose 2}\cdot 4=84$ possibilities to
select such a configuration. It is easy to see that both groups
$\aut(F)$ and $G=\aut(I(F))$ (of order 168 and 336, respectively)
act transitively on the set $\Omega$ of cardinality 84. The
advantage of the group $G$ is that the corresponding association
scheme is of rank 25, that is on the edge of our computational
possibilities. The remaining details relevant to the precise
description of the resulted DSRG$(84,31,17,12,11)$ are in
Appendix.

It is worthy to notice that the incidence graph $I(F)$ is isomorphic to the Heawood graph
considered above. In fact, the considerations from Heawood graph were fulfilled in advance
(this took a few days of computational time) and it was exceptionally extended up to
groups of rank 30, which lead to the discovery of DSRG$(84,29,19,6,12)$ and $(84,39,27,18,18)$.
\end{example}

\begin{example}
Consider the unique $(9,3,1)$-design. We identify its set of points with the set
$\mathcal P=\{1,2,\ldots,9\}$ and its set of blocks with
$$\mathcal B=\{\{1,2,3\},\{1,4,7\},\{1,5,9\},\{1,6,8\},\{2,4,9\},\{2,5,8\},\{2,6,7\},$$
$$\{3,4,8\},\{3,5,7\},\{3,6,9\},\{4,5,6\},\{7,8,9\}\}.$$
The full group of automorphisms of $(\mathcal P,\mathcal B)$ can
be identified with the permutation group $G=\langle
(1,7,3,2,6,9,4,5), (4,6,5)(7,8,9)\rangle$ of order $432$ and
degree 9. (Of course, it could be also regarded as a permutation
group of degree $|\mathcal P|+|\mathcal B|=21$.) Let us consider
the action of $G$ on the 3-set $\{1,2,4\}$. There are $72=12\cdot
3\cdot 2$ possibilities (from geometric arguments) for such a
selection. Denote by $\mathcal O$ the entire set of selected
configurations. The automorphism group of $(\mathcal P,\mathcal
B)$ acts naturally on $\mathcal O$ as a permutation group $\tilde
G$ of degree 72, rank 16. By choosing suitable subsets of 2-orbits
of $\tilde G$, we get two non-equivalent DSRGs with parameter set
$(72,22,9,6,7)$. For a representation of $\tilde G$ as a
permutation group of degree 72 see the group called $H_7$ in
Appendix.
\end{example}

\section{New Sporadic Examples as Cayley Digraphs}

In this section, we construct some new DSRGs of order $32$ and
$39$ as Cayley digraphs. Among them, we obtain the first DSRG with
parameter set $(39,16,12,7,6)$.

The following lemma is crucial for testing whether a Cayley digraph is DSRG.

\begin{lemma}[\cite{km},\cite{hs}]
The Cayley digraph $\cay(G,X)$ is DSRG with parameter set 
\newline 
$(v,k,t,\lambda,\mu)$ if and
only if the equation $\ul{X}\cdot\ul{X} = t\cdot \ul{e} + \lambda\cdot
\ul{X} + \mu\cdot(\ul{G} -\ul{e}-\ul{X})$ holds in~$\mathbb ZG$.
\label{lem1}
\end{lemma}

\subsection{Cayley Digraphs on 32 Vertices}

We now show how to obtain new DSRGs for parameter sets $(32,9,6,1,3)$, $(32,13,9,4,6)$
and $(32,14,10,6,6)$.

Let us take the wreath product group $H=S_2\wr \mathbb Z_4$ of order $32$.
(Here for wreath product we follow notation from \cite{kpr}, which is inherited from
L.A. Kalu\v znin.)
Each element $h\in H$ can be uniquely represented as $h=(g; k_1,k_2)$, where $g\in S_2$, $k_1,k_2\in
\mathbb Z_4$. Let $x$ be a generator of $\mathbb Z_4$, and $\pi=(1\,2)\in S_2$. In order to shorten
description we display just the triple $i\,j\ l$ instead of $(\pi^i;x^j,x^l)$.

Let us define six subsets of $H$:
\begin{eqnarray*}
X_1 &=& \{ 002,011,012,032,033,100,101,102,103\}, \\
X_2 &=& \{ 010,011,030,031,033,100,102,111,113,121,123,130,132\}, \\
X_3 &=& \{ 010,021,022,023,030,100,102,111,113,121,123,130,132\}, \\
X_4 &=& \{ 002,011,012,032,033,100,101,102,103,120,121,122,123\}, \\
X_5 &=& \{ 001,003,011,012,032,033,100,101,102,103,110,111,112,113\}, \\
X_6 &=& \{ 001,003,011,012,032,033,100,101,102,103,110,112,131,133\}.
\end{eqnarray*}

It is a routine-work to check using Lemma \ref{lem1} that the following proposition holds:

\begin{proposition}
The Cayley digraph $\Gamma_i=\cay(H,X_i)$ is a DSRG with parameter set
\begin{itemize}
\item[a)] $(32,9,6,1,3)$, for $i=1$;
\item[b)] $(32,13,9,4,6)$, for $i=2,3,4$, and
\item[c)] $(32,14,10,6,6)$, for $i=5,6$.
\end{itemize}
\label{prop32}
\end{proposition}

According to \cite{ab}, DSRGs with parameter sets mentioned in the previous proposition
had been constructed in \cite{gh}, \cite{ma}.

Using computer algebra system {\sf GAP} \cite{ga} along with the computer package {\sf GRAPE} \cite{gr}
in common with {\sf nauty} \cite{na} we have tested that all the digraphs constructed in Proposition \ref{prop32}
are pairwise non-equivalent and none of these DSRGs were obtained earlier.

\subsection{Cayley Digraphs on 39 Vertices}

In this subsection, we construct DSRGs for all feasible parameter
sets on $39$ vertices. Hence, we obtain also a DSRG for the
parameter set $(39,16,12,7,6)$ for which such a graph has not been
known at the time of writing this paper. All these graphs arise as
Cayley digraphs over a metacyclic group of order $39$. For more
constructions from metacyclic groups, we refer the reader to
\cite{di}.

Let us take the group $G$ presented as
\[
G=\langle a,b: a^{13}=b^3=e, ba = a^9b\rangle\leq AGL(1,13)
\]

and let us define eight of its subsets:

\begin{eqnarray*}
X_1 &=& \{ a,a^5,a^8,a^{12},b,a^2b,a^4b,a^3b^2,a^7b^2,a^{11}b^2\}, \\
X_2 &=& \{ a,a^5,a^8,a^{12},b,a^4b,a^7b,a^{11}b,a^2b^2,a^4b^2,a^8b^2,a^{11}b^2\}, \\
X_3 &=& \{ a,a^5,a^8,a^{12},b,a^4b,a^7b,a^{10}b,a^3b^2,a^4b^2,a^{10}b^2,a^{11}b^2\}, \\
X_4 &=& \{ a^2,a^4,a^9,a^{11},b,a^4b,a^6b,a^{11}b,a^3b^2,a^5b^2,a^9b^2,a^{11}b^2\}, \\
X_5 &=&  X_2\cup\{ a^{12}b, a^3b^2\}, \\
X_6 &=&  X_4\cup\{ a^2b, a^7b^2 \}, \\
X_7 &=& \{ a^2,a^4,a^9,a^{11},a^2b,a^3b,a^6b,a^8b,a^{11}b,a^{12}b, b^2, a^3b^2,a^4b^2,a^6b^2,a^7b^2,a^{10}b^2\}, \\
X_8 &=& \{ a^4,a^5,a^8,a^9,b,ab,a^4b,a^6b,a^8b,a^{10}b,b^2,a^2b^2,a^4b^2,a^6b^2,a^8b^2,a^{10}b^2\}.
\end{eqnarray*}

\begin{proposition}
The Cayley digraph $\Gamma_i=\cay(G,X_i)$ is a DSRG with parameter set
\begin{itemize}
\item[a)] $(39,10,6,1,3)$, for $i=1$;
\item[b)] $(39,12,4,3,4)$, for $i=2,3,4$,
\item[c)] $(39,14,6,5,5)$, for $i=5,6$,  and
\item[d)] $(39,16,12,7,6)$ for $i=7,8$.
\end{itemize}
\label{prop39}
\end{proposition}

\begin{remark}
The DSRG with parameter set $(39,10,6,1,3)$ is isomorphic to the one constructed in \cite{ma}
and described using partial sum families.
The two digraphs with parameters $(39,16,12,7,6)$ are non-equivalent.
\end{remark}

\begin{remark}
In \cite{jo} the author states that it may happen that the so-called Krein parameter
$q_{\theta\theta}^{\theta}$ is always non-negative if $0$ and $-1$ are not eigenvalues
of a DSRG $\Gamma$ in consideration.
The existence of a DSRG with parameters $(39,10,6,1,3)$ disproves it,
since in this case $q_{\theta\theta}^{\theta}=-3/8$ (this fact was somehow not
mentioned in \cite{ma}).
\end{remark}

\section{Conclusion and Summary}

The main genre of this paper is computer algebra experimentation for the purposes
of algebraic graph theory. Using techniques and ideas, which were before reflected
in \cite{di}, \cite{fp}, \cite{fm}, \cite{km} and \cite{kp}, the author \v S. Gy\"urki
arranged a more systematical search for DSRGs, relying on the above described strategies.

We think that the approaches outlined above carry features of
methodological innovations, though in a few cases they simply stem
from careful analysis of previous successful computations done by
M. Klin et al.

Our next goal was of a definite ``sporting'' interest: to present examples of new DSRGs for
previously open parameter sets. Altogether we reached such a success for 28 new parameter sets,
see Table \ref{tabS} below.

Of course, the foremost goal at a computer algebra experimentation
(cf. \cite{zv}) is to reach a successful theoretical
generalization of the obtained new results. In the case of the
sporadic examples of new DSRGs, this would mean to try to embed at
least some of the new examples into new infinite classes of DSRGs.
We are pleased to claim that this task was successfully fulfilled
in the course of our project. In fact, we succeeded to
generalize the presented digraph with parameter set
$(32,14,10,6,6)$ to the infinite series of DSRGs with parameters
$(2n^2,4n-2,2n+2,n+2,6)$. The corresponding paper is in
preparation. Hence, we can finally claim that one more corollary
of the reported project is creation of new (striking in the eyes
of the authors) patterns of successful insight:
\begin{itemize}
\item to observe a short sequence of parameter sets with similar properties;
\item to formulate a plausible conjecture about a possible putative infinite series of
combinatorial structures;
\item to prove this conjecture on purely theoretical level, that is finally, without the use
of a computer.
\end{itemize}

Table~\ref{tabS} below provides a brief summary of our computer
aided discoveries.

\begin{table}[!h]
\caption{Summary. (Used abbreviations: ps -- Is parameter set new?;
am -- The amount of new constructed DSRGs.)
}
\label{tabS}
\begin{center}
\begin{tabular}{|c|c|c|c||c|c|c|c|}
\hline
$(n,k,t,\lambda,\mu)$ & ps & am & constructed in &
$(n,k,t,\lambda,\mu)$ & ps & am & constructed in \\
\hline
(30,13,11,6,5) & Yes & 1 & Section 4.1 & (54,19,9,6,7) & Yes & 2 & Section 4.1\\
(32,9,6,1,3) & No & 1 & Section 5.1 & (54,20,16,6,8) & Yes & 2 & Section 4.1 \\
(32,13,9,4,6) & No & 3 & Section 5.1 & (54,21,17,8,8) & Yes & 2 & Section 4.1 \\
(32,14,10,6,6) & No & 2 & Section 5.1 & (54,25,14,11,12) & Yes & 8 & Section 4.1 \\
(36,13,7,4,5) & Yes & 6 & Section 4.1 & (60,13,5,2,3) & Yes & 1 & Section 4.2 \\
(36,13,11,2,6) & Yes & 1 & Section 4.1 & (60,26,20,10,12) & Yes & 2 & Section 4.3\\
(39,10,6,1,3) & No & 0 & Section 5.2 & (63,22,10,7,8) & Yes & 1 & Section 4.3 \\
(39,12,4,3,4) & No & 3 & Section 5.2 & (72,19,11,2,6) & Yes & 1 & Section 4.3 \\
(39,14,6,5,5) & No & 2 & Section 5.2 & (72,20,14,4,6) & Yes & 1 & Section 4.3 \\
(39,16,12,7,6) & Yes & 2 & Section 5.2 & (72,21,15,6,6) & Yes & 1 & Section 4.3 \\
(45,16,8,5,6) & Yes & 2 & Section 4.1 & (72,22,9,6,7) & Yes & 2 & Section 4.3 \\
(48,10,6,2,2) & Yes & 1 & Section 4.2 & (72,26,10,8,10) & Yes & 1 & Section 4.3 \\
(48,13,7,2,4) & Yes & 1 & Section 4.2 &  (84,29,19,6,12) & Yes & 1 & Section 4.2\\
(50,16,10,3,6) & Yes & 1 & Section 4.1 & (84,31,17,12,11) & Yes & 1 & Section 4.2 \\
(50,23,13,10,11) & Yes & 1 & Section 4.1 & (84,39,27,18,18) & Yes & 2 & Section 4.2 \\
(54,8,3,2,1) & Yes & 1 & Section 4.1 & (90,28,16,10,8) & Yes & 1 & Section 4.3 \\
(54,16,12,6,4) & Yes & 1 & Section 4.1 & (105,36,16,11,13) & Yes & 1 & Section 4.2 \\
\hline
\end{tabular}
\end{center}
\end{table}

\section*{Acknowledgements}

\indent The first author gratefully acknowledges the contribution of the Scientific
Grant Agency of the Slovak Republic under the grant 1/1005/12.

This research was also supported by the Project: Mobility - enhancing research,
science and education at the Matej Bel University, ITMS code: 26110230082,
under the Operational Program Education cofinanced by the European Social Fund.

We thank L. J\o rgensen for generous sharing with us of his
preliminary results related to the DSRG on 108 vertices. A
long-standing cooperation with Ch. Pech and S. Reichard in the use
of computer algebra tools is appreciated.

\section*{Appendix}
\subsection*{Parameter Sets $(48,10,6,2,2)$ and $(48,13,7,2,4)$}
\[
\begin{array}{rl}
H_1=&\langle(1,2,4)(3,11,13,6,7,14)(5,29,22,16,10,26)(8,18,20,9,17,21)(12,24,35)\\
&(15,33,41,37,28,46)(19,32,44)(23,25,36,40,27,38)(30,42,43)(31,34,45,47,39,48),\\
&(1,3)(2,5)(4,8)(6,12)(7,15)(9,19)(10,23)(11,18)(13,25)(14,26)(16,30)\\
&(17,31)(20,33)(22,34)(24,36)(27,39)(32,45)(35,46)(37,47)(38,44)(41,42)(43,48)\rangle
\end{array}
\]
Digraph with vertex set $\{1,2,\ldots,48\}$ and dart set 
$(1,3)^{H_1}\cup(1,5)^{H_1}\cup(1,15)^{H_1}$ is a DSRG $(48,10,6,2,2)$, while the one 
with dart set 
$(1,2)^{H_1}\cup(1,5)^{H_1}\cup(1,12)^{H_1}\cup(1,19)^{H_1}\cup(1,23)^{H_1}$ is
a DSRG $(48,13,7,2,4)$.

\subsection*{Parameter Set $(60,13,5,2,3)$}
\[
\begin{array}{rl}
H_2=&\langle(1,5,16,40,21)(2,8,24,14,37)(3,10,29,28,9)(4,12,26,39,15)\\
&(6,18,44,33,54)(7,20,36,49,23)(11,31,53,58,60)(13,35,34,42,32) \\
&(17,41,27,50,43)(19,45,51,57,25)(22,46,55,56,48)(30,52,38,59,47), \\
&(1,23,40,3,43,34)(2,15,32,6,28,12)(4,37,22,11,57,8)(5,35,46,14,41,58)\\
&(7,21,49,19,47,56)(9,54,30,25,60,17)(10,16,53,27,24,51)\\
&(13,39,52,33,20,55)(18,42,50,31,48,59)(26,29,45,36,44,38)\rangle.
\end{array}
\]
Digraph with vertex set $\{1,2,\ldots,60\}$ and dart set 
$(1,5)^{H_2}\cup(1,6)^{H_2}\cup(1,7)^{H_2}\cup(1,8)^{H_2}\cup(1,12)^{H_2}
\cup(1,19)^{H_2}\cup(1,20)^{H_2}\cup(1,49)^{H_2}$
is a DSRG with parameter set $(60,13,5,2,3)$. 

\subsection*{Parameter Set $(60,26,20,10,12)$}
\[
\begin{array}{rl}
H_3=&\langle(1,23,54,6,49,40)(2,32,50,14,57,37)(3,31,27,12,45,43)(4,16,41,9,30,28)\\
&(5,22,46,24,26,53)(7,11,58,21,39,47)(8,10,36,17,19,60)(13,15,20,33,35,52)\\
&(18,48,59,44,38,42)(25,29,51,34,56,55),\\
&(1,56,31,34,16,55)(2,48,15,18,32,42)(3,24,60,46,41,47)(4,38,10,44,23,59)\\
&(5,33,53,20,58,37)(6,39,45,21,30,26)(7,29,22,25,11,51)(8,17,43,27,54,40)\\
&(9,57,19,14,49,35)(12,52,36,50,28,13)\rangle.\\
\end{array}
\]
Digraphs with vertex set $\{1,2,\ldots,60\}$ and dart set 
$(1,3)^{H_3}\cup(1,5)^{H_3}\cup(1,12)^{H_3}\cup(1,17)^{H_3}\cup(1,36)^{H_3}$,
$(1,12)^{H_3}\cup(1,17)^{H_3}\cup(1,23)^{H_3}\cup(1,36)^{H_3}$
are two non-equivalent DSRGs with parameter set $(60,26,20,10,12)$.

\subsection*{Parameter Set $(63,22,10,7,8)$}
\[
\begin{array}{rl}
H_4=&\langle(1,3,8,19,38,33,61)(2,5,12,7,17,34,20,29,50,16,32,25,45,62)\\
&(4,9,21,11,18,36,6,15,30,24,42,55,13,27)(10,23,41,26,35,40,28,48,63,53,37,56,54,44)\\
&(14,39,57,52,43,31,22)(46,47,59,49,60,51,58),\\
&(1,2)(3,6)(5,10,50,61)(7,8)(9,14)(11,24,42,51)(12,19)(16,35)(17,46)\\
&(18,31)(20,30)(21,58)(23,39,63,29)(25,55)(26,32)(27,41)(28,33)(36,47)(37,43)(48,57)\\
&(52,56)(53,59)(60,62).
\end{array}
\]
Digraph with vertex set $\{1,2,\ldots,63\}$ and dart set 
$(1,2)^{H_4}\cup(1,3)^{H_4}\cup(1,5)^{H_4}$ is a DSRG with 
parameter set $(63,22,10,7,8)$.

\subsection*{Parameter Set $(72,19,11,2,6)$}
\[
\begin{array}{rl}
H_5=&\langle(1,6,47,35,8,23,56,24)(2,13,36,32,4,40,64,41)(3,10,61,46,19,15,54,16)\\
&(5,18,17,65,12,26,30,57)(7,21,51,44,11,28,63,29)(9,31,25,72,20,43,42,68)\\
&(14,53,49,66,27,50,33,59)(22,38,60,71,39,34,45,70)(48,52)(55,58,69,67),\\ 
&(1,28,42,37,70,44,8,10,45,55,72,41)(2,15,30,52,71,46,4,21,33,67,68,24)\\
&(3,40,25,48,59,32,19,6,60,58,65,29)(5,34,64,12,38,63,20,53,54,9,50,56)\\
&(7,23,17,62,66,35,11,13,49,69,57,16)(14,26,36,22,43,47,39,31,61,27,18,51), \\
&(1,15,51,72,48,12)(2,28,61,65,37,5)(3,23,36,68,62,20)(4,13,64,59,69,39)\\
&(6,56,70,67,22,8)(7,40,47,57,52,9)(10,54,71,55,14,19)(11,21,63,66,58,27)\\
&(16,50,60,35,18,30)(17,24,34,49,32,26)(25,41,38,33,44,43)(29,53,45,46,31,42)\rangle.\\
\end{array}
\]
Digraph with vertex set $\{1,2,\ldots,72\}$ and dart set 
$(1,2)^{H_5}\cup(1,7)^{H_5}\cup(1,10)^{H_5}\cup(1,42)^{H_5}$
is a DSRG with parameter set $(72,19,11,2,6)$.

\subsection*{Parameter Sets $(72,20,14,4,6)$ and $(72,21,15,6,6)$}
\[
\begin{array}{rl}
H_6=&\langle(1,61,53,14,31,33,63,48)(2,41,26,16,59,42,51,20)(3,71,28,21,34,38,49,13)\\
 &(4,39,44,40,35,6,52,36)(5,64,11,9,69,22,50,30)(7,25,60,68,17,10,58,23) \\
 &(8,72,19,32,47,29,56,27)(12,43,70,57,24,18,67,15)(37,45,65,55,62,46,54,66), \\
&(2,18,3)(4,24,37)(5,10,8)(7,27,45)(9,20,54)(11,33,25)(12,13,46)(14,17,62)\\
&(15,29,31)(16,30,65)(19,22,34)(21,48,55)(23,38,39)(26,44,43)(28,42,47)(32,35,66)\\
&(36,70,59)(40,51,56)(41,58,61)(49,63,50)(52,64,67)(53,60,69)\rangle,\\
\end{array}
\]
Digraph with vertex set $\{1,2,\ldots,72\}$ and dart set 
$(1,7)^{H_6}\cup(1,21)^{H_6}\cup(1,40)^{H_6}\cup(1,57)^{H_6}$
is a DSRG with parameter set $(72,20,14,4,6)$.

Digraph with vertex set $\{1,2,\ldots,72\}$ and dart set 
$(1,2)^{H_6}\cup(1,4)^{H_6}\cup(1,15)^{H_6}\cup(1,36)^{H_6}$
is a DSRG with parameter set $(72,21,15,6,6)$.

\subsection*{Parameter Set $(72,22,9,6,7)$}
\[
\begin{array}{rl}
H_7=&\langle(1,44,55,3,20,25,48,28)(2,36,64,42,19,51,69,49)(4,37,9,34,13,39,17,38) \\
&(5,53,47,68,14,26,66,71)(6,18,59,10,12,31,40,16)(7,46,29,57,8,56,70,62)\\
&(11,72,67,15,58,23,45,54)(21,41,30,50,61,43,33,27)(22,35,65,24,52,60,32,63),\\
&(1,32,51)(2,43,53)(3,45,30)(4,69,31)(5,52,19)(6,33,60)(7,61,20)(8,16,72)\\
&(9,48,66)(10,68,65)(11,56,26)(12,15,50)(13,17,70)(14,49,58)(18,47,25)(21,71,37)\\
&(22,62,46)(23,39,36)(24,34,35)(27,67,38)(28,63,54)(29,40,44)(41,42,57)(55,59,64)\rangle,
\\
\end{array}
\]

Digraphs with vertex set $\{1,2,\ldots,72\}$ and dart set 
$(1,2)^{H_7}\cup(1,10)^{H_7}\cup(1,14)^{H_7}\cup(1,24)^{H_7}\cup(1,71)^{H_7}$, 
and $(1,2)^{H_7}\cup(1,10)^{H_7}\cup(1,24)^{H_7}\cup(1,45)^{H_7}\cup(1,72)^{H_7}$
are two non-equivalent DSRGs with parameter set $(72,22,9,6,7)$.

\subsection*{Parameter Set $(72,26,10,8,10)$}
\[
\begin{array}{rl}
H_8=&\langle(1,10,9,12,25,18)(2,49,50,46,41,59)(3,63,11,47,20,69)(4,31,33,60,52,66)\\
&(5,56,71,51,38,53)(6,36,17)(7,35,57,19,27,39)(8,24,54,61,14,40)(15,45,68)\\
&(13,26,43,28,67,58)(16,23,55)(21,65,30,48,37,32)(22,44,42,34,62,72)(29,64,70), \\
&(1,58)(2,67)(3,26)(4,68)(5,45)(6,72)(7,13)(8,15)(9,53)(10,59)(11,28) \\
&(12,70)(14,49)(16,43)(17,30)(18,35)(19,38)(20,51)(21,55)(22,65)(23,44) \\
&(24,66)(25,56)(27,47)(29,54)(31,48)(32,60)(33,39)(34,57)(36,63)(37,41) \\
&(40,62)(42,61)(46,71)(50,69)(52,64), \\
&(1,63,36,58,5,45)(2,40,18,59,48,13)(3,54,6,56,21,68)(4,55,25,72,29,26)\\
&(7,31,10,35,62,67)(8,23,9,30,64,28)(11,52,17,53,44,15)(12,70,51,43,16,20)\\
&(14,57,22,37,39,66)(19,42,69,27,60,71)(24,33,41,65,34,49)(32,47,50,61,38,46)\rangle.\\
\end{array}
\]
Digraph with vertex set $\{1,2,\ldots,72\}$ and dart set 
$(1,3)^{H_8}\cup(1,14)^{H_8}\cup(1,19)^{H_8}\cup(1,20)^{H_8}\cup(1,58)^{H_8}$
is a DSRG with parameter set $(72,26,10,8,10)$.

\subsection*{Parameter Sets $(84,29,19,6,12)$, $(84,31,17,12,11)$ and $(84,39,27,18,18)$}

\[
\begin{array}{rl}
H_9=&\langle
 (1,42,52,14,67,78,38,39)(2,28,54,40,18,37,41,23)(3,64,77,82,80,53,33,56)\\
&(4,34,30,15,79,83,61,58)(5,63,74,27,32,72,68,31)(6,13,36,12)\\
&(7,11,25,29,10,35,43,45)(8,51,75,69,50,55,26,49)(9,20,57,84,81,70,62,22)\\
&(16,24,46,76,71,73,17,21)(19,44,59,48,47,60,65,66), \\
&(1,71,8)(2,72,20)(3,10,60)(4,68,23)(5,82,14)(6,63,52)(7,75,24)(9,25,28)\\
&(11,16,54)(12,44,74)(13,42,49)(15,67,77)(17,73,39)(18,78,30)(19,81,27)\\
&(21,38,40)(22,65,29)(26,59,36)(31,58,80)(32,84,57)(33,69,45)(34,41,37)\\
&(35,56,43)(46,51,76)(47,79,64)(48,55,83)(50,53,61)(62,70,66)\rangle.
\end{array}
\]

Digraph with dart set $(1,3)^{H_9}\cup(1,5)^{H_9}\cup
(1,6)^{H_9}\cup(1,7)^{H_9}\cup(1,8)^{H_9}\cup(1,9)^{H_9}\cup
(1,19)^{H_9}\cup(1,26)^{H_9}\cup(1,29)^{H_9}\cup(1,40)^{H_9}$
is a DSRG$(84,29,19,6,12)$,
digraph with dart set $(1,2)^{H_9}\cup(1,3)^{H_9}\cup
(1,5)^{H_9}\cup(1,7)^{H_9}\cup(1,19)^{H_9}\cup(1,20)^{H_9}\cup
(1,27)^{H_9}\cup(1,30)^{H_9}\cup(1,31)^{H_9}\cup(1,35)^{H_9}$
is a DSRG$(84,31,17,12,11)$, 
digraphs with dart set $(1,3)^{H_9}\cup(1,4)^{H_9}\cup
(1,5)^{H_9}\cup(1,6)^{H_9}\cup(1,13)^{H_9}\cup(1,14)^{H_9}\cup
(1,19)^{H_9}\cup(1,26)^{H_9}\cup(1,29)^{H_9}\cup(1,30)^{H_9}
\cup(1,35)^{H_9}\cup(1,40)^{H_9}\cup(1,46)^{H_9}$
and
$(1,3)^{H_9}\cup(1,4)^{H_9}\cup
(1,5)^{H_9}\cup(1,6)^{H_9}\cup(1,13)^{H_9}\cup(1,15)^{H_9}\cup
(1,19)^{H_9}\cup(1,21)^{H_9}\cup(1,23)^{H_9}\cup(1,26)^{H_9}
\cup(1,29)^{H_9}\cup(1,30)^{H_9}\cup(1,40)^{H_9}$
are DSRG$(84,39,27,18,18)$.

\subsection*{Parameter Set $(90,28,16,10,8)$}
\[
\begin{array}{rl}
H_{10}=&\langle(1,52,41,43,40)(2,15,79,34,24)(3,75,64,66,56)(4,63,60,14,31)(5,86,69,76,36)\\
&(6,20,13,82,83)(7,25,73,48,12)(8,39,88,37,45)(9,33,32,62,74)(10,70,51,59,17)\\
&(11,80,53,65,67)(16,61,58,90,55)(18,26,72,21,27)(19,81,28,42,49)\\
&(22,71,78,54,68)(23,47,46,85,44)(29,50,38,87,77)(30,84,57,89,35)\\
&(1,79,84,2,60,78)(3,73,61,7,83,77)(4,41,53)(5,34,50,10,48,68)\\
&(6,70,47,12,49,38)(8,51,23,15,72,36)(9,25,88,17,18,69)(11,26,75,20,40,55)\\
&(13,64,28)(14,86,33,24,67,22)(16,82,80,27,74,85)(19,39,52,31,56,35)(21,32,58)\\
&(29,43,71,42,59,89)(30,63,81,45,44,62)(37,46,57)(54,66,87,65,76,90)\rangle.\\
\end{array}
\]
Digraph with vertex set $\{1,2,\ldots,90\}$ and dart set 
$(1,2)^{H_{10}}\cup(1,5)^{H_{10}}\cup(1,15)^{H_{10}}\cup(1,20)^{H_{10}}\cup(1,40)^{H_{10}}$
is a DSRG with parameter set $(90,28,16,10,8)$.

\subsection*{Parameter Set $(105,36,16,11,13)$}
\[
\begin{array}{rl}
H_{11}=&\langle(1,96,45,103,3,68,12,105,32,79)(2,82,50,100,10,42,15,102,53,60)(13,25,77,97,49)\\
&(4,47,28,104,36,56,5,85,64,98)(6,99,22,93,7,87,20,101,11,90)(16,39,88,83,44)\\
&(8,91,38,84,17,69,24,95,29,74)(9,73,40,92,14,70,18,94,26,81)(21,72,89,65,27)\\
&(19,55,57,86,43,41,30,78,71,67)(23,34,61,75,48)(31,46,62,80,37)\\
&(33,35,66,76,63)(51,54,52,58,59),\\
&(1,97)(2,89)(3,88)(4,77)(5,83)(6,76)(7,80)(8,62)(9,63)(10,72)(11,64)(12,65)\\
&(13,66)(14,50)(15,49)(16,52)(17,59)(18,58)(19,104)(20,53)(21,51)(22,54)(23,102)\\
&(24,36)(25,37)(26,35)(27,46)(28,44)(29,45)(30,100)(31,38)(32,40)(33,39)(34,98)\\
&(41,105)(42,95)(43,94)(47,92)(48,91)(55,103)(56,85)(57,84)(60,82)(61,81)(67,99)\\
&(68,101)(69,73)(70,74)(71,75)(78,93)(79,96)\rangle
\end{array}
\]

Digraph with vertex set $\{1,\ldots,105\}$ and dart set $(1,2)^{H_{11}}\cup(1,4)^{H_{11}}
\cup(1,10)^{H_{11}}\cup(1,22)^{H_{11}}\cup(1,26)^{H_{11}}$ is a DSRG with parameter set $(105,36,16,11,13)$.

\end{document}